\documentclass[reqno,12pt]{amsart}
\usepackage{graphicx,a4wide}
\usepackage{mathrsfs}

\numberwithin{equation}{section}

\newtheorem{theo}{Theorem}

\newtheorem{coro}{Corollary}

\newtheorem{lem}{Lemma}
\theoremstyle{remark}

\def\({\left(}
\def\){\right)}
\def\[{\left[}
\def\]{\right]}

\def\dis{\displaystyle}
\def\dd{\textup{d}}

\newcommand\bin[2]{\left( \begin{array}{c}#1\\ #2\end{array}\right)}
\newcommand\beq{\begin{equation}}
\newcommand\eeq{\end{equation}}

\date{}

\begin{document}

\title[]{A family of criteria for irrationality of Euler's constant }
\author[]{Marc Pr\'evost  }
  
\email{prevost@lmpa.univ-littoral.fr}

\address{LMPA Joseph Liouville, Centre Universitaire de la Mi-Voix,
B\^at. H. Poincar\'e, 50 rue F. Buisson, BP 699, 62228 Calais Cedex, France} 
 
 

\begin{abstract}
	Following earlier results of Sondow, we propose another criterion of irrationality for   Euler's constant $\gamma$. It involves similar linear combinations of logarithm numbers $L_{n,m}$. To prove that $\gamma$ is irrational, it suffices to prove that, for some fixed $m$, the distance of  $d_n L_{n,m}$   ($d_n$ is the least common multiple of the $n$ first integers )  to the set of integers $\mathbf{Z}$ does not converge to $0$.
	A similar result is obtained by replacing logarithms numbers by  rational numbers: it gives  a sufficient condition involving only rational numbers.
	Unfortunately, the chaotic behavior of $d_n$  is an obstacle to verify this sufficient condition.\\
	All the proofs use in a large manner the theory of  Pad\'e approximation.
\end{abstract}

\maketitle
\section{Introduction }
 
In \cite{sondow1}, the author, using Beukers' integral \cite{beukers}, found a criterion for irrationality of Euler's constant $\gamma$. It depends on the limit of the fractional part of the following expression
$$ L_n=2\sum_{k=1}^n \sum_{i=0}^{k-1} \bin{n}{i}^2(H_{n-i}-H_i)\ln(n+k)$$
where $H_n$ is the Harmonic number $\dis H_n:=\sum_{k=1}^n \frac{1}{k}$.
In this paper, we establish the connexion between Sondow's criteria and Pad\'e approximant of the function $\dis \frac{\ln {u}}{u-1}
$.
Moreover, following the same idea, we find a family of new criteria:
\emph{for each integer $m \leq n$, let us  set
 \beq\label{definition:Lnm}
L_{n,m}:=\sum_{k=0}^n\bin{n}{k}\bin{n+k}{k}(-1)^{n+k}\ln(n-m+k+1),
\eeq
$d_n:={\bf LCM}(1,\ldots,n) $ and  $\{x\} =x- \left\lfloor x\right\rfloor $  the fractional part of the real $x$. 
If, for some integer $m$, if the sequence $\{d_n (-1)^m L_{n,m}\}$ does not converge to $0$ when $n$ tends to infinity, then $\gamma$ is irrational.}

Using the property of the error term, a more precise criterion is proved here:\\
if, {\it for some integer $m$,   the sequence $\{d_{2^n} (-1)^m L_{2^n,m}\}$ is asymptotically non  decreasing, when $n$ tends to infinity, then $\gamma$ is irrational.}

\section{Sondow's criterion with Pad\'e approximant}\label{sondow:criterion}
Sondow considers the double integral (so-called Beukers' integral)
$$I_n=\int_0^1\int_0^1\frac{(x\,(1-x)\,y\,(1-y))^n}{(1-x\, y)\ln (\,x\, y)}\; \dd x \; \dd y.$$

Applying Taylor expansion of $1/(1-xy)$ around $0$, he proved the following identity

\begin{eqnarray} \label{sondowformula}
	I_n=\bin{2n}{n} \gamma +L_n-A_n=\mathcal{O}(2^{-4n}n^{-1/2})
\end{eqnarray}

where  
$ \dis A_n=\sum_{i=0}^n\bin{n}{i}^2H_{n+i}$.
After multiplication by $d_{2n}$, it arises
$$d_{2n}\bin{2n}{n}\gamma=d_{2n}(A_n-L_n)+o(1).$$
{\bf Sondow's criterion:}

Since $d_{2n}A_n \in {\bf Z}$,  if the sequence of  fractional part $(\{d_{2n} L_n\})_n$ does not converge to $0$ then  $\gamma \notin \bf Q$.
 
\vskip 1cm

Sebah computed this sequence for $1\leq n\leq 2500$. Its cumulative average   seems to converge $1/2$, but the mathematical proof remains to establish.

In the following, we will show that the sequence involved in the paper by Sondow can be recovered by means of Pad\'e approximation.

Let us consider the function $(\ln u)/(u-1)$ and its Pad\'e Approximant  $[n-1/n]$ of degree $(n-1/n)$ at the point $u=1$:
\beq \label{formule:padelog}
\frac{\ln u	}{u-1}=\frac{N_n(u)}{D_n(u)}+R_n(u)
\eeq
where $N_n$ and $D_n$ are polynomials of respective degree $n-1$ and $n$, normalized by $N_n(1)=D_n(1)=1$, and $R_n(u)=\mathcal{O}(u^{2n})$

From the theory of Pad\'e approximation, it is well known  that $D_n$ is related with the shifted Legendre Polynomial orthogonal on the interval [0,1] with respect to the Lebesgue weight function. Some of these      expressions 	are
\begin{eqnarray}
P_n^*(t)&=&\sum_{k=0}^n \bin{n}{k}^2 t^{n-k}(t-1)^k\label{formulePn:1}\\
&=&\sum_{k=0}^n \bin{n}{k} \bin{n+k}{k} (-1)^{n+k} t^k \label{formulePn:2}
\end{eqnarray}
$D_n$ has the following expression in terms of $P_n^*$
\begin{equation}\label{formule:1}
	D_n(u)=P^*_n\(\frac{1}{1-u}\)(1-u)^n\bin{2n}{n}^{-1}.
\end{equation}

Replacing  $P_n^*$ by its expressions (\ref{formulePn:1},\ref{formulePn:2}), formula (\ref{formule:1}) becomes
\begin{eqnarray*}
		D_n(u)&=&\bin{2n}{n}^{-1}\sum_{k=0}^n\bin{n}{k}^2 u^k\\ &=&\bin{2n}{n}^{-1}\sum_{k=0}^n\bin{n}{k} \bin{n+k}{k} (1-u)^{n-k}u^k.
\end{eqnarray*}

 The numerator $N_n(u)$ of  $[n-1/n]$, is related with the associated polynomial of the denominator:
\begin{eqnarray*}
	N_n(u)&=&2 \bin{2n}{n}^{-1}\sum_{k=1}^n \sum_{i=0}^{k-1} \bin{n}{i}^2
(H_{n-i}-H_i)u^{k-1}\\
	&=&\bin{2n}{n}^{-1}\sum_{k=0}^n\bin{n}{k}\bin{n+k}{k}\sum_{i=0}^{k-1}(u-1)^{n-k+i} \frac{(-1)^i}{i+1} .
\end{eqnarray*}
\\
\\

Now, what is the link between    Pad\'e Approximation and Sondow's criterion?\\
\\
The definition of $\gamma$ is, primarily, $$\gamma=\lim_n(H_n-\ln n).$$
An integral representation for Euler's constant  is 
\begin{equation}\label{formule1:gamma}
	\gamma=\int_0^1\(\frac{1}{\ln u}+\frac{1}{1-u}\) \; \dd u .
\end{equation}

Formula (\ref{formule:padelog}) can be rewritten as
\begin{equation*}
\frac{N_n(u)}{\ln u}+\frac{D_n(u)}{1-u}=-\frac{R_n(u)D_n(u)}{\ln u}
\end{equation*}
and 
\begin{equation*}
\gamma=\int_0^1\frac{1-u^n N_n(u)}{	\ln u} \;\dd u +\int_0^1\frac{1-u^n D_n(u)}{	1-u} \;\dd u
-\int_0^1\frac{u^n R_n(u)D_n(u)}{\ln u}\;\dd u
\end{equation*}
By linearity, the second term is expanded as
\begin{eqnarray} 
\int_0^1\frac{1-u^n D_n(u)}{	1-u} \;\dd u  &= &\bin{2n}{n}^{-1}\sum_{k=0}^n\bin{n}{k}^2 \int_0^1
\frac{u^{n+k}-1}{u-1}\;\dd u \\
&=&\bin{2n}{n}^{-1}\sum_{k=0}^n\bin{n}{k}^2 
H_{n+k} \\
&=&\bin{2n}{n}^{-1} A_n.\label{formule:An}
\end{eqnarray}

 The first integral can be computed as following:
 \begin{eqnarray}
\int_0^1\frac{1-u^n N_n(u)}{\ln(u)} \dd u &
= &2\bin{2n}{n}^{-1}\sum_{k=1}^n \sum_{i=0}^{k-1} \bin{n}{i}^2
(H_{n-i}-H_i) \int_0^1 \frac{1-u^{n+k-1}}{\ln(u)}  \dd u\\
&=&-2 \bin{2n}{n}^{-1}\sum_{k=1}^n \sum_{i=0}^{k-1} \bin{n}{i}^2
(H_{n-i}-H_i)  \ln(n+k)\\
&=&-\bin{2n}{n}^{-1} L_n.\label{formule:Ln}
\end{eqnarray}

 From the theory of Pad\'e approximation, in formula (\ref{formule:padelog}), the remainder term  $R_n$ has an integral representation
   $$R_n(u)=\frac{(1-u)^n}{D_n(u)}\int_0^1 \frac{t^nD_n(1-1/t)}{1-(1-u)t}\; \dd t.$$
   Thanks to formulas (\ref{formule:An},\ref{formule:Ln}), $\gamma$ satisfies
\begin{eqnarray*}
\gamma&=&\frac{A_n-L_n}{\bin{2n}{n}}-\int_0^1  \frac{u^n R_n(u)D_n(u)}{\ln(u)}\;\dd u.\\
\end{eqnarray*}

Thus another expression of the remainder term $I_n$ of Sondow is

\begin{eqnarray*}
I_n&=&\bin{2n}{n}\gamma-A_n+L_n =-\int_0^1  \frac{u^n R_n(u)D_n(u)}{\ln(u)}\; \dd u\\
&=&\int_0^1 \int_0^1 \frac{u^n(1-u)^n}{\ln(u)}\frac{P_n^*(t)}{1-(1-u)t}\; \dd t \; \dd u\\
&=&\int_0^1 \int_0^1 \frac{u^n(1-u)^{2n}}{\ln(u)}\frac{t^n(1-t)^n}{(1-(1-u)t)^{n+1}}\; \dd t \; \dd u
\end{eqnarray*}
thanks to integration by parts and Rodrigues formula for orthogonal polynomials.

Thus the approximation for Euler's constant $\gamma$ (\ref{sondowformula}) is a consequence of the Pad\'e approximation to the function $\dis (\ln u)/(1-u)$.

In the same manner, Pilehrood \cite{pilehrood} found irrationality criteria for generalized Euler's constant. He defined the following linear form in logarithms \begin{eqnarray*}
L_{(n_1,n_2)}(\alpha)&=&\sum_{m=1}^{n_1} \sum_{k=0}^{m-1} \bin{n_1}{k} \bin{n_2}{k}(H_{n_1-k}+H_{n_2-k}-2H_k)\ln(m+n_1+\alpha-1)+\\&&
\sum_{m=n_1+1}^{n_2} \sum_{k=m}^{n_2} (-1)^{k-1-n_1}/k\bin{n_2}{k}/ \bin{k-1}{n_1}(\ln(m+n_1+\alpha-1)
\end{eqnarray*}

Actually, following the same idea as for Sondow's criterion, it is possible to prove that Pilehrood' criterion comes from Pad\'e approximations $[n_2-1,n_1]=\frac{R_{n_2-1}(u)}{S_{n_1}(u)} $ (normalized by $ R_{n_2-1}(1)=S_{n_1}(1)=1$ to the function $\ln(u)/(u-1) $ at the point $u=1$. The linear form $L_{(n_1,n_2)}(\alpha)$ satisfies:
$$L_{(n_1,n_2)}(\alpha)=\int_0^1\left(1-u^{n_1+\alpha-1}{R_{n_2-1}(u)} \right)\frac{1}{\ln(u)}\dd u$$

\section{Statement of the results}
 In order to simplify Sondow's criterion, it is convenient to choose a more simple approximation.
 This method leads to the following theorems.

 \begin{theo}\label{theo:0}
 For $0\leq m\leq n, $ let us define 
 
\begin{eqnarray*}
	J_{n,m}&=&\int_0^1 u^{n-m} P_n^*(u)\dis \left( \frac{1}{\ln u}+\frac{1}{1-u}\right)\; \dd  u\\
L_{n,m}&=&-\int_0^1 \(\frac{1-u^{n-m}P_n^*(u)}{\ln u}\)\; \dd  u\\
A_{n,m}&=&\int_0^1 \(\frac{1-u^{n-m}P_n^*(u)}{1- u}\)\; \dd  u
\end{eqnarray*}
then 
\begin{eqnarray}
\gamma&=&A_{n,m}-L_{n,m}+J_{n,m}\label{formuletheo1}\\
 	L_{n,m}
	&=& \sum_{k=0}^n \bin{n}{k} \bin{n+k}{k} (-1)^{n+k}   \ln(n-m+k+1)\\
	A_{n,m}&=&2H_n
\end{eqnarray}

 \end{theo}

 \begin{theo}\label{theo:1}
 
 The following are equivalent:\\
 (a) The fractional part of $d_n L_{n,m}$ is given by $\{d_n (-1)^m L_{n,m}\}=d_n (-1)^m J_{n,m} \;\; (*)$ for some $n$ or $m$.\\
 (b) The formula (*) holds for some $m$ and  for all sufficiently large $n$.\\
 (c) Euler's constant is a rational number.\\
 \end{theo}

 A sufficient condition which   involves $d_n L_{n,m}$ but not  $J_{n,m}$ 
 is the following
 \begin{theo}\label{theo:2}
 If for some integer $m$, $\{d_n (-1)^m L_{n,m}\}\geq 0.707^n $ infinitely often, then $\gamma$ is irrational.
 \end{theo}
 
 Computations (see Table \ref{table:1}) show that this condition is satisfied for  $n\leq 1000$.
 Numerical results also suggest   that for each $m$,  $\{d_n (-1)^m L_{n,m}\}$ is dense in the interval (0,1) and the cumulative average $n^{-1}\sum_{k=1}^\infty \{d_k (-1)^m L_{k,m}\}$converges to $0.5 $ (see Figure 1 and 2).
 
 To prove $\gamma$ irrational, it just  suffices   to show that $\{d_n (-1)^m L_{n,m}\}$ does not converge to $0$.

\begin{figure}
	\centering		\includegraphics[width=0.80\textwidth]{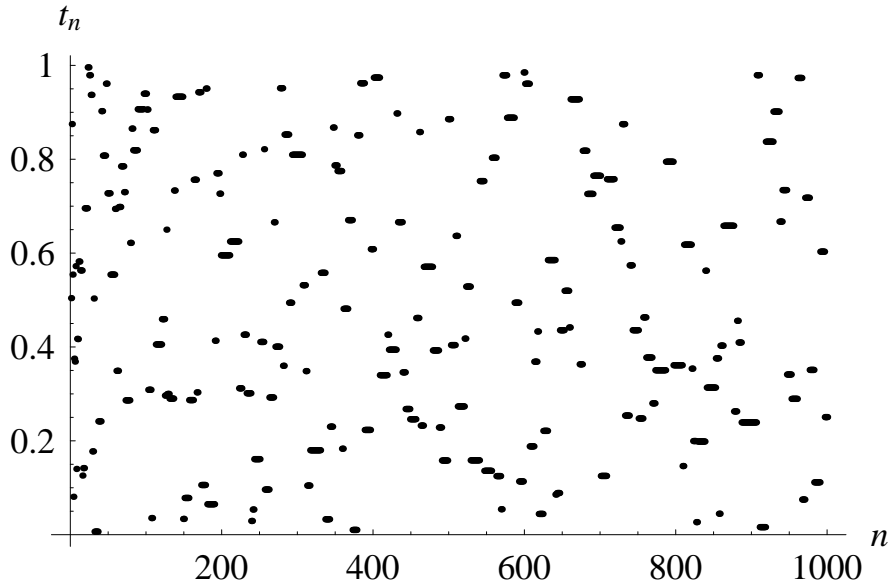}
		\caption{plot of values of $t_n=\{d_n (-1)^m L_{n,m}\}$, for $m=0$ }
	\label{fig:figure2}
\end{figure}

\begin{figure*}
	\centering
		\includegraphics[width=0.80\textwidth]{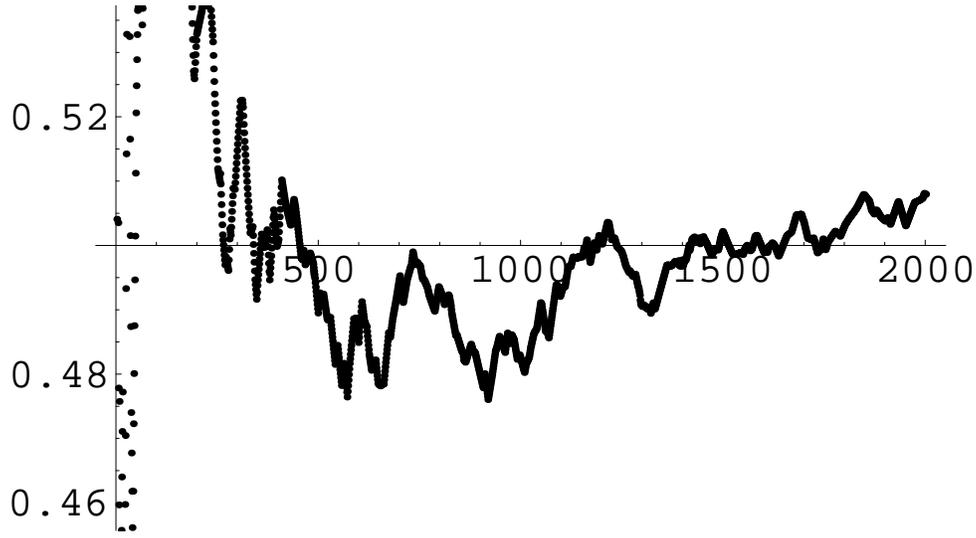}
			\caption{plot of values of $<t_n>= n^{-1}\sum_{k=1}^n\; t_k$}
	\label{fig:moyenne}
\end{figure*}

  \begin{table}  \caption{}
\begin{tabular}{|l  l  l l l|}
\hline
&&&&\\
$ n$ &  $\dis \frac{.7^{n}}{\left\{ d_{n} L_{n,0}\right\} }$ &$\dis \frac{.7^{n}}{\left\{ - d_{n}L_{n,1}\right\} }$  &$\dis \frac{.7^{n}}{\left\{ d_{n}L_{n,2}\right\} }$&$\dis \frac{.7^{n}}{\left\{ - d_{n}L_{n,3}\right\} }$\\ \hline
&&&&\\
 1& 1.38868 & 1.81209 & --------- &---------\\
2 & 0.56003 & 0.58439 & 0.56609 & ---------\\
3 & 0.61882 & 0.64252 & 0.63428 &0.67030 \\
4 & 2.97160 & 3.31151 & 3.23310 &0.38225 \\
5 & 0.44808 & 0.45886 & 0.45719 &0.45913 \\
6 & 0.31896 & 0.32064 & 0.32044 &0.32061 \\
7 & 0.14391 & 0.14467 & 0.14460 &0.14465 \\
8 & 0.41138 & 0.41543 & 0.41511 &0.41528 \\
9 & 0.09667 & 0.09689 & 0.09687 &0.09688 \\
10& 0.06778 & 0.06781 & 0.06781 &0.06781 \\
11& 0.03395 & 0.03398 & 0.03398 &0.03398 \\
12& 0.02378 & 0.02379 & 0.02379 &0.02379 \\
13& 0.01719 & 0.01721 & 0.01720 &0.01720 \\
14& 0.01204 & 0.01204 & 0.01204 &0.01204 \\
15& 0.00843 & 0.00843 & 0.00843 &0.00843 \\
16& 0.02637 & 0.02637 & 0.02637 &0.02637 \\
17& 0.01639 & 0.01639 & 0.01639 &0.01639 \\
18& 0.01147 & 0.01147 & 0.01147 &0.01147 \\
19& 0.00163 & 0.00163 & 0.00163 &0.00163 \\
20&0.001147 & 0.00114 & 0.00114 &0.00114\\
  \hline
\end{tabular}\label{table:1}
\end{table}

In section \ref{section:6}, we will prove the asymptotic formula
\beq\label{asymptformula}
\gamma=A_{n,m}-L_{n,m}+\mathcal{O}(4^{-n}).
\eeq
Actually, we will prove that the error term   $J_{n,m}$  is 
  a totally monotone sequence (i.e. a sequence of moments with respect a positive measure), converging to $0$ as $4^{-n}$.\\
  
 By substituting in $L_{n,m}$, $\ln (n+1+k-m)$ by some suitable Pad\'e approximants, a sufficient condition, involving only rational numbers is the following
 
 \begin{coro}\label{coro:1}
  Let us define, for $n-m+1=2^p,p\in {\mathbf Z}$,
 \beq\label{definition:Ltildenm}
 \tilde{L}_{n,m}:=   p \;[n/n]_{t=1}+\sum_{k=0}^n\bin{n}{k}\bin{n+k}{k}(-1)^{n+k}[n/n]_{ t=k/(n-m+1)}
\eeq
where $[n/n] $ is the Pad\'e approximant of $\ln(1+t)$ at $t=0$:
\begin{equation*}
[n/n]_{t}=
\dis
  \frac{\dis t
      \sum_{k=0}^n \bin{n}{ k} 
         \bin{n+k}{k}  \left(
         \sum_{i=0}^{k-1} \frac{t^{i-k+n}(-1)^i}{i+1} \right)}
         {
         \dis
      \sum_{k=0}^n \bin{n}{ k} 
         \bin{n+k}{k}t^{n-k}  }  
\end{equation*}

 If for some integer $m$, $\{d_{2^p} (-1)^m \tilde{L}_{2^p+m-1,m}\}$ does not converge to $0$ when $p$ tends to infinity, then $\gamma$ is irrational.

 \end{coro}
 
 Another sufficient condition comes from the property of the error term in the asymptotic formula (\ref{asymptformula}) and from the upper and lower bound of the $LCM(1,\ldots,n)$:
 \begin{coro}\label{coro:2}
 
  If for some $m$, $\{d_{2^p} (-1)^m {L}_{2^p,m}\}$ is asymptotically non decreasing, then $\gamma$ is irrational.
 \end{coro}
 
\section{Two lemmas}

\begin{lem}\label{lemma:1}

 The function $\dis \frac{1}{\ln(1- u)}+\frac{1}{u}$ is a Markov-Stieltjes function.
More precisely, 
\beq\label{logMS}
 \frac{1}{\ln (1-u)}+\frac{1}{u}=\int_0^1 \frac{1}{1-u\; t}w(t) \; \dd  t
\eeq
  where the weight function $w$ is 
 $$w(t):=\frac{1}{t\;(\ln^2(1/t-1)+\pi^2)}$$

 \end{lem}
 
 \begin{proof}

 After a change of variable ($u\rightarrow (1-u)$ and $x=1/t-1$), formula (\ref{logMS}) is equivalent to 
 
\beq\label{logMS2}
 \frac{1}{\ln (u)}+\frac{1}{1-u}=\int_0^\infty \frac{1}{x+u}\frac{1}{\ln^2x+\pi^2} \; \dd t
\eeq
 The weight function $w$ can be found with the Stieltjes inversion formula (see \cite{widder}). Another  way to prove formula (\ref{logMS2}) is to apply residue theorem to the function $$f(x):=\frac{1}{x+u}\frac{1}{\ln x+i \pi}.$$
 Taking the determination of $\ln x$ on the complex plane cut along the positive real axis, the poles of $f$ are $x=-u$ and $x=-1$.
 
 Let us define $\gamma_r$ a small semi-circle $z=r e^{i \theta}, -\pi/2\leq\theta\leq \pi/2, r>0$.
 $D^+_r$ the line $z=x+i r, x  {\rm \;running \;from }\;0 \; to \; R $, $\Gamma_R$ the circle $z=R e^{i \theta},0\leq \theta\leq 2 \pi$ and $D^-_r$ the line $z=x-i r, {\rm  \;for }\;x  {\rm  \;from }\;R \; to \; 0 $.

 Now, we  compute $\int_{\mathcal C}f(x)\; \dd x$ where ${\mathcal C}$ is the union of $D^+_r$, $\Gamma_R$, $D^-_r$ and $\gamma_r$, with the theorem of residue to obtain
 \begin{eqnarray*}
 \int_0^\infty
\frac{1}{x+u}\(\frac{1}
{\ln x+i\pi}+\frac{-1}{\ln x-i\pi}\)\; \dd x&= &\int_0^\infty
\frac{1}{x+u}\(\frac{-2i\pi}{\ln^2 x+\pi^2 }\)\; \dd x\\
&=&-2 i \pi
\(\frac{1}{\ln u}+\frac{1}{1-u}\)
 \end{eqnarray*}
 \end{proof}
 
 \vskip 1cm
 \noindent Now, we are in position to prove a new formula for the Euler's constant $\gamma$.

 \begin{theo}
 
  The Euler's constant $\gamma$ satisfies
 \begin{equation*}
 \gamma=\int_{-\infty}^{\infty}\frac{\ln (1+e^{-z})e^z}{z^2+\pi^2}\; \dd  z
\end{equation*}
\end{theo}

 \begin{proof}
 
  In the integral representation of $\gamma$ (\ref{formule1:gamma}), let us substitute 
 the integrand by the expression (\ref{logMS}).
  This leads to 
  \begin{eqnarray*} 
  \gamma&=&\int_0^1-\frac{\ln (1-t)}{t}\frac{1}{t(\ln^2(1/t-1)+\pi^2)}\;\; \dd t\\
  &=&\int_{-\infty}^{\infty}\frac{\ln (1+e^{-z})e^z}{z^2+\pi^2}\;\dd z
\end{eqnarray*}
with the change of variable $t=(1+e^z)^{-1}$.
 
 \end{proof}

\begin{lem}\label{lemma:2}
For each fixed integer $m$, the sequence $((-1)^mJ_{n,m})_n$ defined in Theorem  \ref{theo:1} is  totally monotonic. More precisely
\begin{equation*}
(-1)^m J_{n,m}= \int_0^{1/4} v^n \rho_m(v)\; \dd v
\end{equation*} 
where the weight function is $$ \rho_m(v)=\int_{\frac{1-\sqrt{1-4\;v}}{2}}^{\frac{1+\sqrt{1-4\;v}}{2}}\(\frac{u-u^2-v}{u\;v}\)^m\frac{1}{(u-u^2-v)\(\dis \pi^2+\ln^2\(\frac{-u\;v}{u^2-u+v}\)\)}\;\dd u.$$
\end{lem}

 \begin{proof}
  $\dis J_{n,m}=\int_0^1 u^{n-m} P_n^*(u) \(\frac{1}{\ln u}+\frac{1}{1-u}\)\; \dd  u$
  appears as Legendre modified moments of the weight function $\dis \(\frac{1}{\ln u}+\frac{1}{1-u}\)\; \dd  u$.
  \\
For some particular cases of weight function, a sequence of modified moments can be itself a sequence of moments, with respect to a positive measure (see \cite{prevost2}).
Using Rodrigues formula for orthogonal polynomials, Lemma \ref{lemma:1}, Fubini's theorem and after $n$ integrations by parts, it arises 
\begin{eqnarray*}
	J_{n,m}&=&\int_0^1 u^{n-m} \frac{(-1)^n}{n!}\frac{\dd ^n}{\dd u^n}\(u^n(1-u)^n\) \(\frac{1}{\ln u}+\frac{1}{1-u}\)\; \dd  u\\
	&=&\int_0^1 u^{n-m} \frac{(-1)^n}{n!}\frac{\dd ^n}{\dd u^n}\(u^n(1-u)^n\) \;\dd u
	\int_0^1 \frac{1}{1-(1-u)\; t}w(t) \;\dd t \\
	&=&\int_0^1 \int_0^1  \frac{(-1)^n}{n!}u^n(1-u)^n \;\dd u
	 \frac{\dd ^n}{\dd u^n}\(\frac{u^{n-m}}{1-(1-u)\; t}\)\; w(t) \;\dd t .
\end{eqnarray*}
  
  The computation of $\dis \frac{\dd ^n}{\dd u^n}\(\frac{u^{n-m}}{1-(1-u)\; t}\)$ needs the partial decomposition of the rational function $\dis \frac{u^{n-m}}{1-(1-u)\; t}=q(u)+\dis \(\frac{t-1}{t}\)^{n-m}\;\frac{1}{1-(1-u)t}$,
  where $q$ is polynomial of degree $n-m-1$.

Another expression of $J_{n,m}$ is then
\begin{eqnarray*}
	J_{n,m}&=&\int_0^1 \int_0^1   u^n(1-u)^n 
	   \(\frac{t-1}{t}\)^{n-m}\frac{t^n}{(1-(1-u)t)^{n+1}}\; w(t) \;\dd t \; \dd u.
\end{eqnarray*}
  We do the following change of variable $$\dis v=\frac {u(1-u)(1-t)}{1-(1-u)t} \in [0,1/4]
  \Leftrightarrow t=\phi(v)=\frac{u^2-u+v}{(u-v)(u-1)} \in[0,1].$$
  Let $\phi_1(v)$ and $\phi_2(v)$ denote the two roots of the quadratic equation
  $v=u-u^2$,\\ $\dis  \phi_1(v)~=~\frac{1+\sqrt{1-4\;v}}{2},$ $\dis \phi_2(v)=\frac{1-\sqrt{1-4\;v}}{2}$.

\begin{eqnarray*}
	J_{n,m} 
	&=& \int_0^{1/4} v^n\; \dd v\int_{\phi_1(v)}^{\phi_2(v)}(-1)^m \(\frac{\phi(v)	}{\phi(v)-1}\)^m w(\phi(v)) \frac{u^2}{(u-1)(u-v)^2}\;\frac{(-1)^m}{1-(1-u)\phi(v)}\;\dd u
\end{eqnarray*}
which proves the lemma.

\end{proof}

\section{Proof of Theorem~\ref{theo:0}} 
We first prove the  identity (\ref{formuletheo1}) linking Euler's constant $\gamma$, the linear combination of logarithms numbers $L_{n,m}$, the rational numbers $A_{n,m}$ and the integrals $J_{n,m}$.
From formula (\ref{formule1:gamma}), one substitute the integrand $\dis \(\frac{1}{\ln u}+\frac{1}{1-u}\)$ by an approximation involving Legendre Polynomials as follows:
\begin{eqnarray*}
\gamma&=&\int_0^1 \(\frac{1}{\ln u}+\frac{1}{1-u}\)\; \dd  u\\
&=&\int_0^1 \(\frac{1-u^{n-m}P_n^*(u)}{\ln u}+\frac{1-u^{n-m}P_n^*(u)}{1-u}\)\; \dd  u+\\
&&\int_0^1 u^{n-m} P_n^*(u)\dis \left( \frac{1}{\ln u}+\frac{1}{1-u}\right)\; \dd  u.
\end{eqnarray*}
The   expression (\ref{formulePn:2}) of $P_n^*$ leads to analogous expressions $L_{n,m}$.\\
 By linearity
\begin{eqnarray*}
	L_{n,m}&= &-\sum_{k=0}^n \bin{n}{k} \bin{n+k}{k} (-1)^{n+k} \int_0^1 \(\frac{1-u^{k+n-m} }{\ln u}\)\\
	&=& \sum_{k=0}^n \bin{n}{k} \bin{n+k}{k} (-1)^{n+k}   \ln(n-m+k+1)
\end{eqnarray*}

\noindent $A_{n,m}$  
 is treated quite differently: \\
 $P_n^*$ satisfies the following orthogonality relation
 
\begin{eqnarray*}
	\int_0^1 P_n^*(u)q(u)\; \dd  u=0, \rm{for \; all}  \rm{ \; polynomial \;} \emph{{q}} \rm{ \; of\; degree \;less \;than  \;} \emph{n}.
\end{eqnarray*}
Thus, by taking $\dis q(u)=\frac{1-u^{n-m}}{1-u}$,   another expression for 	$A_{n,m}$ is 
\begin{eqnarray}\label{formule2:Anm}
		A_{n,m}=\int_0^1\frac{1-P_n^*(u)}{1-u}\; \dd  u=\int_0^1\frac{P_n^*(1)-P_n^*(u)}{1-u}\; \dd  u
\end{eqnarray}
	and so $A_{n,m}$ is independent of $0\leq m\leq n$.
	
	\noindent Let us now compute the integral in (\ref{formule2:Anm}).
	
	\noindent	Legendre polynomials satisfy a three term recurrence relation which is
	
\begin{eqnarray*}
	(n+1)P_{n+1}^*(u)=(2n+1)(2u-1)P_n^*(u)-n P_{n-1}^*(u)\\
	P_0^*(u)=1\;\;\;
	P_1^*(u)=2u-1	\hskip 4cm
\end{eqnarray*}

Thus,   $A_{n,m}$ 's also satisfy a similar recurrence relation

\begin{eqnarray}
	(n+1)A_{n+1,m}=(2n+1)(2u-1)A_{n,m}(u)-n A_{n-1,m}(u)\label{rec1A}\\
	A_0=0\;\;\;
	A_1=2	\hskip 4cm\label{rec2A}
\end{eqnarray}
With (\ref{rec1A}) and (\ref{rec2A}), it is not difficult to prove that 
\begin{eqnarray*}
	A_{n,m}=2 H_n,0\leq m \leq n
\end{eqnarray*}

\section{Proof of Theorem~\ref{theo:1} and Theorem~\ref{theo:2}}\label{section:6}
All the arguments are based on the formula (\ref{formuletheo1}).
\begin{eqnarray*}
	\gamma=A_{n,m}-L_{n,m}+J_{n,m}\Longleftrightarrow d_n 	\gamma=d_n A_{n,m}-d_n L_{n,m}+d_n J_{n,m}
\end{eqnarray*}
Thus, 
\begin{eqnarray}
	d_n 	\gamma \in \mathbf{Z} &\Longleftrightarrow& d_n L_{n,m}-d_n J_{n,m} \in \mathbf{Z}\label{irratcriterion1}\\ 
  &\Longleftrightarrow& \{d_n (-1)^m L_{n,m}   \}=\{d_n (-1)^m J_{n,m}\}\label{irratcriterion2}   \end{eqnarray}
 since $A_{n,m}=2 H_n$ and $(-1)^m J_{n,m} $ is positive.
 
 On the other hand, Lemma \ref{lemma:2} implies that the sequence $\(J_{n,m}\)_n$ converges to $0$ as $4^{-n}$. The numbers $d_n$ converges to infinity as $e^n$.
 Thus $d_n (-1)^m J_{n,m}$ is decreasing  to $0$. So,  for all sufficiently large $n$,    $\{d_n (-1)^m J_{n,m}\} =d_n (-1)^m J_{n,m}$.
 
 If $\gamma$ is a rational number, then $ d_N \gamma \in \mathbf{Z}$ for some $N$ and  for all $n\geq N$, $d_n \gamma \in \mathbf{Z}$. $(c)\Longrightarrow (b)$  follows from the previous arguments.
 Thus $(c)\Longrightarrow (b) \Longrightarrow (a) \Longrightarrow (c)$.

In  (\ref{irratcriterion2}), we substitute $d_n$ by an upper bound : $d_n \leq e^{1.039\; n}$ \cite{rosser}.\\ Thus $\forall n, \;d_n (-1)^m J_{n,m} \leq e^{1.039 \; n}\; 4^{-n} < 0.707^n$ and Theorem~\ref{theo:2} is proved.

\section{Proof of Corollaries}

\noindent 1) In the numerical computation of  formula $$L_{n,m}= \sum_{k=0}^n \bin{n}{k} \bin{n+k}{k} (-1)^{n+k}   \ln(n-m+k+1),$$ the problem is the evaluation of logarithmic functions.\\
A mean to avoid this drawback is the substitution of $\ln(n-m+k+1)$ by some suitable approximations, enough good to keep the irrationality criteria.
We will show now that Pad\'e approximants satisfy this condition:\\
another expression of $L_{n,m}$ is 
\begin{equation*}
	L_{n,m}= \ln (n-m+1)+\sum_{k=0}^n \bin{n}{k} \bin{n+k}{k} (-1)^{n+k}   \ln\(1+\frac{k}{n-m+1}\).
\end{equation*}
The Pad\'e error for the logarithmic function is 
\begin{eqnarray}\label{padeerror}
	\ln (1+x)-[n/n]_x=\frac{(-1)^n x^{n+1}}{P_n^*(-1/x)}\int_0^1\frac{t^n\;(1-t)^n}{(1+x\;t)^{n+1}}\;\dd t
\end{eqnarray}
   Let us set 
   
\begin{eqnarray*}
	 {L'}_{n,m}:=   \ln (n-m+1)+\sum_{k=0}^n\bin{n}{k}\bin{n+k}{k}(-1)^{n+k}[n/n]_{ t=k/(n-m+1)} 
\end{eqnarray*}

 We have to evaluate the difference   $\delta_{n,m}:=L_{n,m}-{L'}_{n,m}$.
 For sake of simplicity, we set $\dis \zeta_k=\frac{k}{n-m+1}$.
 
\begin{eqnarray*}
	\delta_{n,m}  &  =   & \sum_{k=0}^n \bin{n}{k} \bin{n+k}{k} (-1)^{n+k} \(  \ln\(1+\zeta_k\) - [n/n]_{ t=\zeta_k} \)\\
	&=&\sum_{k=0}^n \bin{n}{k} \bin{n+k}{k} (-1)^{k} 
 \frac{\zeta_k^{n+1}}{P_n^*(-\zeta_k^{-1})}
\int_0^1\frac{t^n\;(1-t)^n}{(1+\zeta_k\;t)^{n+1}}\;\dd t	\\
 \end{eqnarray*}

 \noindent Since $\zeta_k \in [0,1]$ and $P_n^*$ has all its roots in $[0,1]$, $\dis \left|\frac{\zeta_k^{n}}{P_n^*(-\zeta_k^{-1})}\right|\leq \frac{1}{\left|P_n^*(-1)\right|}$.
 On the other hand, the integral  
\begin{eqnarray*}
	\int_0^1\frac{t^n\;(1-t)^n}{(1+\zeta_k\;t)^{n+1}}&\leq&  4^{-n}\int_0^1 \frac{1} {(1+\zeta_k\;t)^{n+1}}\\	 
	&\leq&  4^{-n} \frac{1}{n\; \zeta_k}.
\end{eqnarray*}
 So, 
 \begin{eqnarray*}
\left|\delta_{n,m}\right| & \leq &\sum_{k=0}^n \bin{n}{k} \bin{n+k}{k}  
 \left|\frac{\zeta_k^{n+1}}{P_n^*(-\zeta_k^{-1})}\right|  \left|
\int_0^1\frac{t^n\;(1-t)^n}{(1+\zeta_k\;t)^{n+1}}\;\dd t  \right| 
 \\& \leq &
\sum_{k=0}^n \bin{n}{k} \bin{n+k}{k}  
 \frac{\zeta_k}{\left|P_n^*(-1)\right|}   4^{-n} \frac{1}{n\; \zeta_k} \\
 & \leq &\frac{1}{\left|n P_n^*(-1)\right|}\;4^{-n}\sum_{k=0}^n \bin{n}{k} \bin{n+k}{k}  = \frac{1}{\left|n P_n^*(-1)\right|}\;4^{-n}\left|P_n^*(-1)\right|=(n \; 4^n)^{-1}
 \end{eqnarray*}
 
 The goal is partly reached since the error between $L_{n,m}$ and its approximation is less than $J_{n,m}$.
   Now, let us consider the approximation of $\ln (n-m+1)$. 
   It is difficult to approximate this number (which tends to infinity) with an error lees than $4^{-n}$. So, we consider sequences of integers $n$, such that $n-m+1$ is a power of $2$: $n-m+1=2^p$. With this hypothesis, $\ln (n-m+1)=p\; \ln 2$
   
    In (\ref{padeerror}), if $x=1$, 
    	$\dis \ln 2-[n/n]_{x=1}=\frac{(-1)^n}{P_n^*(-1)}\int_0^1\frac{t^n\;(1-t)^n}{(1+\;t)^{n+1}}\;\dd   t$.
    	The asymptotics for Legendre polynomials are well known $$ P_n(\alpha)\sim (\alpha+\sqrt{\alpha^2-1})^n, {\rm \; for \;} \alpha \in \mathbf{R}\setminus [-1,1].$$ Thus shifted Legendre Polynomials satisfy

    	$$ P_n^*(t)\sim ((2 t-1)+2\sqrt{ t^2-t})^n, {\rm \; for \;} t \in \mathbf{R}\setminus [0,1].$$

    	    	The maximum of the fraction $\dis \(\frac{t(1-t)}{1+t}\) $ for $t\in [0,1]$ is  
    	obtained for $t=\sqrt{2}-1$, and its  value is $(3-2 \sqrt{2})$. Thus
\begin{eqnarray*}
\left|	\ln 2-[n/n]_{x=1}\right| \leq  \frac{(3-2 \sqrt{2})^n}{(3+2\sqrt{2})^n}\; \ln 2
\end{eqnarray*}

 \noindent  For $n-m+1=2^p$, $\ln (2^p)-p \; [n/n]_{x=1} \leq p\;(3-2\sqrt{2})^{2\; n} $ which is a $ o(4^{-n}/n)$.
At last, the error $\left|L_{n,m}-\tilde{L}_{n,m}\right|$ satisfies
$$\left|L_{n,m}-\tilde{L}_{n,m}\right|\leq (4^{-n}/n)$$
and the corollary 1 is proved.

\noindent 2) For the proof of Corollary \ref{coro:2}, we exploit the property of totally monotonic sequences (TMS). \\
A sequence $u_n$ is called TMS if there exists a non negative measure  $ d \mu $ with infinitely  many points of increase such that $$\forall n \in \mathbf{N},\; u_n=\int_0^\infty x^n \; \dd\mu(x).$$
If the support of the measure $d \mu$ is the interval $[0,1/R]$, then $ \forall n,\dis u_{n+1}/u_n\leq R$ and 
$\dis \lim_n \frac{u_{n+1}}{u_n}=R$. 
If $R=1$, it  is equivalent to $$\forall n \in \mathbf{N},  \forall k \in \mathbf{N}, (-1)^k\Delta^k(u_n)>0$$
where $\Delta^0(u_n):=u_n$ and $\Delta^{k+1} u_n=\Delta^{k}u_{n+1}-\Delta^{k}u_n$.
 (see \cite{widder}, p. 108).
 
 The previous properties can be applied to the sequence $J_{n,m}$ for which we prove some convergence properties. If they are not satisfied by $\{d_n (-1)^m L_{n,m}\}$ then $\gamma$ is irrational.
 
 First we will prove that $J_{n,m}$ satisfies $d_{2n} (-1)^m J_{2n,m}< d_{n} (-1)^m J_{n,m}$:
the numbers $d_n$ and $J_{n,m}$ satisfy $2^n \leq d_n<e^{1.039 \;n}$ (see \cite{tenenbaum}, p.12-13 for the lower bound and \cite{rosser} for the upper one) $\dis \frac{J_{n+1,m}}{J_{n,m}}<1/4$ (property of totally monotonic sequence \cite{widder}, p.135).

 $$\frac{d_n  J_{n,m}}{d_{2n}J_{2n,m}}
 >\frac{2^n}{e^{1.039 \times \;2n}} \; {4^n}> 1.0014  $$
 Thus, for all integer $m$, $\(  d_{2^p}(-1)^m J_{2^p,m}\)_{p\in \mathbf{N}}$ is a positive decreasing sequence, converging to $0$. So, if 
 $\(\{d_{2^p} (-1)^m L_{2^p,m}\}\)_p$ is non decreasing for $p$ greater than any integer,   then $\gamma$ is irrational.

 \vskip 1cm
 I would thank my colleague S. Eliahou for reference  (\cite{pilehrood})

\def\refname{Bibliography}


\begin{thebibliography}{1}
\bibitem[Beukers 1979]{beukers}F. Beukers,
\textit{A note on the irrationality of $\zeta (2)$ and $\zeta (3)$},
Bull. London Math. Soc. \textbf{11}, 268-272 (1979).
\bibitem[Rosser et al. 1962]{rosser} J. Rosser, L. Schoenfeld,
 \textit{Approximate formulas for some functions of prime numbers},
Illinois J. Math. \textbf{6}, 64-94 (1962).   
\bibitem[Sondow 2003]{sondow1}{J. Sondow, \textit{Criteria for irrationality of Euler's constant }, Proc. Amer. Math. Soc. \textbf{131} (2003), 3335-3344.}
\bibitem[Pilehrood 2004]{pilehrood} T. Hessami Pilehrood, Kh. Hessami Pilehrood, \textit{Criteria for irrationality of generalized Euler's constant}, J. of Number Theory, 108 (2004) 169-185.
\bibitem[Pr\'evost 1994]{prevost2} M. Pr\'evost, 
\textit{Acceleration of some logarithmic sequences}, 
J. Comput. Appl. Math. \textbf{55}, No.3, 357-367 (1994).
\bibitem[Tenenbaum 1990]{tenenbaum} G. Tenenbaum,
\textit{Introduction \`a la th\'eorie analytique et probabiliste des nombres}, 
 Institut Elie Cartan, Universit\'e de Nancy   (1990).
\bibitem[Widder 1941]{widder} D. V. Widder, 
\textit{The Laplace transform},
Princeton Mathematical Series,  Princeton University Press, N. J. X, 406 p. (1941). 
  
\end{thebibliography}
\end{document}